\documentclass[12pt]{article}
\usepackage{hyperref}
\usepackage{amssymb, amsmath}
\usepackage{graphicx}
\usepackage[mathscr]{eucal}

\begin{document}
\title{Why there is no an existence theorem\\
       for a convex polytope with prescribed\\
       directions and perimeters of the faces?}
\author{Victor Alexandrov}
\date{Originally submitted: July 25, 2017. Corrected: November 25, 2017}
\maketitle
\begin{abstract}
We choose some special unit vectors 
$\boldsymbol{n}_1,\dots,\boldsymbol{n}_5$
in $\mathbb{R}^3$ and denote by $\mathscr{L}\subset\mathbb{R}^5$ 
the set of all points $(L_1,\dots,L_5)\in\mathbb{R}^5$ 
with the following property:
there exists a compact convex polytope $P\subset\mathbb{R}^3$ 
such that the vectors $\boldsymbol{n}_1,\dots,\boldsymbol{n}_5$ 
(and no other vector) are unit outward normals to the faces of 
$P$ and the perimeter of the face with the outward normal 
$\boldsymbol{n}_k$ is equal to $L_k$ for all $k=1,\dots,5$.
Our main result reads that $\mathscr{L}$ is not 
a locally-analytic set, i.\,e., we prove that, for some point 
$(L_1,\dots,L_5)\in\mathscr{L}$, 
it is not possible to find a neighborhood $U\subset\mathbb{R}^5$ 
and an analytic set $A\subset\mathbb{R}^5$ such that 
$\mathscr{L}\cap U=A\cap U$. 
We interpret this result as an obstacle for finding an 
existence theorem for a compact convex polytope with prescribed 
directions and perimeters of the faces.
\end{abstract}

\textbf{1. Introduction and the statement of the main result.}
In 1897, Hermann Minkowski proved the following uniqueness 
theorem:
 
\textbf{Theorem 1} (H. Minkowski, \cite{Mi97} and \cite[p. 103--121]{Mi11}).
\textit{A convex polytope is uniquely determined, up to
translations, by the directions and the areas of its faces}.

Here and below a convex polytope is the convex hull of a finite 
number of points.
By  the direction of a face, we mean the direction of the 
outward normal to the face. 

Theorem 1 has numerous applications and generalizations. 
In order to discuss some of them, we will use the following 
notation.

Let $P$ be a compact convex polytope in $\mathbb{R}^3$ and
$\boldsymbol{n}\in\mathbb{R}^3$ be a unit vector. 
By $P^{\boldsymbol{n}}$ we denote the intersection of $P$
and its support plane with the outward normal $\boldsymbol{n}$. 
Note that $P^{\boldsymbol{n}}$ is either a vertex, or an edge, 
or a face of $P$. Accordingly, we say that $P^{\boldsymbol{n}}$ 
has dimension 0, 1, or 2.

In 1937, A.D. Alexandrov proved several generalizations of 
the above 
uniqueness theorem of Minkowski, including the following

\textbf{Theorem 2} (A.D. Alexandrov, \cite{Al37} and \cite[p. 19--29]{Al96}). 
\textit{Let $P_1$ and $P_2$ be convex polytopes in $\mathbb{R}^3$.
Then one of the following mutually exclusive possibilities realizes:}

(i) \textit{$P_1$ is obtained from $P_2$ by a parallel translation};  

(ii) \textit{there exist $k=1,2$ and a unit vector 
$\boldsymbol{n}\in\mathbb{R}^3$ such that
$P_k^{\boldsymbol{n}}$ has dimension 2 and, for some 
translation $T:\mathbb{R}^3\to\mathbb{R}^3$, the formula 
$T(P_j^{\boldsymbol{n}})\subsetneq P_k^{\boldsymbol{n}}$ 
holds true, where $j\in\{1,2\}\diagdown\{k\}$.}

Note that the formula 
$T(P_j^{\boldsymbol{n}})\subsetneq P_k^{\boldsymbol{n}}$ 
means that the face $P_j^{\boldsymbol{n}}$ can be 
embedded inside the face $P_k^{\boldsymbol{n}}$ by translation 
$T$ as a proper subset.

For more details about Theorems 1 and 2, the reader is referred to \cite{Al05}.

For us, it is important that Theorem 1 is a special case of Theorem 2.
In fact, suppose the conditions of Theorem 1 are fulfilled. 
Then, using the notation of Theorem 2, we observe that, 
for every unit vector $\boldsymbol{n}$ such that $P_k^{\boldsymbol{n}}$ 
has dimension 2, $P_j^{\boldsymbol{n}}$ also has dimension 2 
and its area is equal to the area of $P_k^{\boldsymbol{n}}$. 
Hence, there is no translation $T$ such that
$T(P_j^{\boldsymbol{n}})\subsetneq P_k^{\boldsymbol{n}}$. 
This means that the possibility (ii) in Theorem 2
is not realized. Therefore, the possibility (i) in Theorem 2
is realized and Theorem 1 is a consequence of Theorem 2.

In fact, Theorem 2 has many other consequences, including the following 

\textbf{Theorem 3} (A.D. Alexandrov, \cite[Chapter II, \S\;4]{Al05}).
\textit{A convex polytope in $\mathbb{R}^3$ is uniquely 
determined, up to
translations, by the directions and the perimeters of its faces}.

For the sake of completeness, we mention that a direct analog 
of Theorem 1 is valid
in $\mathbb{R}^d$ for all $d\geqslant 4$; a direct analog 
of Theorem 2 is not valid
in $\mathbb{R}^d$ for every $d\geqslant 4$; in $\mathbb{R}^3$, 
a refinement of Theorem 2 
was found by G.Yu. Panina \cite{Pa08} in 2008.

We explained above that uniqueness Theorems 1 and 3 are similar 
to each other and both follow from Theorem 2. 
In the rest part of this section  we explain the difference 
that appears when we are interested in existence results 
corresponding to uniqueness Theorems 1 and 3.

In 1897, Hermann Minkowski also proved the following existence 
theorem:
 
\textbf{Theorem 4} (H. Minkowski, \cite{Mi97} and \cite[p. 103--121]{Mi11}).
\textit{Let unit vectors 
$\boldsymbol{n}_1,\dots,\boldsymbol{n}_m$ in $\mathbb{R}^3$
and real numbers $F_1,\dots,F_m$ satisfy the following conditions:} 

(i) \textit{$\boldsymbol{n}_1,\dots,\boldsymbol{n}_m$ 
are not coplanar and no two of them coincide with each other};

(ii) \textit{$F_k$ is positive for every $k=1,\dots,m$};

(iii) $\sum_{k=1}^m F_k \boldsymbol{n}_k=0$.

\noindent\textit{Then there exists a convex polytope 
$P\subset\mathbb{R}^3$ such that 
$\boldsymbol{n}_1,\dots,\boldsymbol{n}_m$ (and no other vector)
are outward face normals for $P$ and $F_k$ is the area of 
the face with
outward normal $\boldsymbol{n}_k$ for every $k=1,\dots,m$.}

For the sake of completeness, we mention that a direct 
analog of Theorem 4 is valid
in $\mathbb{R}^d$ for all $d\geqslant 4$. 

For more details about Theorem 4, the reader is referred 
to \cite{Al05}.

Recall that a set $A\subset\mathbb{R}^d$ is said to be algebraic 
if $A=\{x\in\mathbb{R}^d: p(x)=0\}$ for some polynomial
$p:\mathbb{R}^d\to\mathbb{R}$, 
and $A$ is said to be locally-algebraic if, for every $x\in A$, 
there is a neighborhood
$U\subset\mathbb{R}^d$ and an algebraic set 
$A_0\subset\mathbb{R}^d$ such that $A\cap U=A_0\cap U$.
 
For a given set $\{\boldsymbol{n}_1,\dots,\boldsymbol{n}_m\}$ 
of vectors satisfying the condition (i) of Theorem 4, 
denote by $\mathscr{F}(\boldsymbol{n}_1,\dots,\boldsymbol{n}_m)$
the set of all points $(F_1,\dots,F_m)\in\mathbb{R}^m$ 
such that there exists a convex polytope $P\subset\mathbb{R}^3$
for which $\boldsymbol{n}_1,\dots,\boldsymbol{n}_m$ (and no 
other vector) are the outward face normals for $P$, and $F_k$ 
is the area of the face with the outward normal 
$\boldsymbol{n}_k$ for every $k=1,\dots,m$.
The set $\mathscr{F}(\boldsymbol{n}_1,\dots,\boldsymbol{n}_m)\subset\mathbb{R}^m$
can be referred to as a natural configuration space of convex 
polytopes (treated up to translations) with prescribed set 
$\{\boldsymbol{n}_1,\dots,\boldsymbol{n}_m\}$ 
of outward unit normals when a polytope is determined by
the areas $F_1,\dots,F_m$ of its faces.

From Theorem 4, it follows immediately that 
\textit{the set $\mathscr{F}(\boldsymbol{n}_1,\dots,\boldsymbol{n}_m)\subset\mathbb{R}^m$
is locally-algebraic for 
every set $\{\boldsymbol{n}_1,\dots,\boldsymbol{n}_m\}$ 
of vectors satisfying the condition} (i) \textit{of Theorem~4.}
In fact, we can define the algebraic set $A_0$ as the zero set
of the quadratic polynomial 
$$
\sum_{j=1}^{3} \biggl(\sum_{k=1}^m 
(\boldsymbol{n}_k,\boldsymbol{e}_j) F_k\biggr)^2,
$$
where $(\boldsymbol{n}_k,\boldsymbol{e}_j)$ stands for the
standard scalar product in $\mathbb{R}^3$ and 
$\{\boldsymbol{e}_1, \boldsymbol{e}_2, \boldsymbol{e}_3\}$
is the standard orthonormal basis in $\mathbb{R}^3$.

Recall that a set $B\subset\mathbb{R}^d$ is said to be analytic 
if $B=\{x\in\mathbb{R}^d: \varphi(x)=0\}$ for some 
real-analytic function $\varphi:\mathbb{R}^d\to\mathbb{R}$, 
and $B$ is said to be locally-analytic if, for every $x\in B$, 
there is a neighborhood $U\subset\mathbb{R}^d$ and an analytic 
set $B_0\subset\mathbb{R}^d$ such that $B\cap U=B_0\cap U$.

Obviously, every locally-algebraic set is locally analytic.

Let unit vectors $\boldsymbol{n}_1,\dots,\boldsymbol{n}_5$
in $\mathbb{R}^3$ be defined by the formulas
\begin{equation}
\begin{split}
\boldsymbol{n}_1=(0,0,-1),\quad
\boldsymbol{n}_2=&\biggl(\frac{1}{\sqrt{2}},0,\frac{1}{\sqrt{2}}\biggr),\quad
\boldsymbol{n}_3=\biggl(-\frac{1}{\sqrt{2}},0,\frac{1}{\sqrt{2}}\biggr),\\
\boldsymbol{n}_4=&\biggl(0,\frac{1}{\sqrt{2}},\frac{1}{\sqrt{2}}\biggr),\quad
\boldsymbol{n}_5=\biggl(0,-\frac{1}{\sqrt{2}},\frac{1}{\sqrt{2}}\biggr).{}
\end{split}
\end{equation}
For convenience of the reader, the vectors $\boldsymbol{n}_1,\dots,\boldsymbol{n}_5$
are shown schematically in Figure~1.
\begin{figure}
\begin{center}
\includegraphics[width=0.35\textwidth]{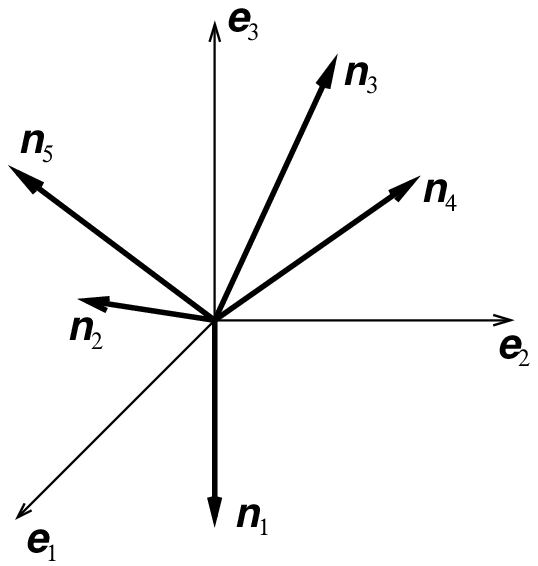}
\end{center}
\caption{Unit vectors $\boldsymbol{n}_1,\dots,\boldsymbol{n}_5$ 
defined by the formulas (1)}
\end{figure}
By $\mathscr{L}=\mathscr{L}(\boldsymbol{n}_1,\dots,\boldsymbol{n}_5)\subset\mathbb{R}^5$ 
we denote the set of all points $(L_1,\dots,L_5)\in\mathbb{R}^5$ 
with the following property:
there exists a convex polytope $P\subset\mathbb{R}^3$ such that 
the vectors $\boldsymbol{n}_1,\dots,\boldsymbol{n}_5$ (and no 
other vector) are the unit outward normals to the faces of $P$, 
and $L_k$ is the perimeter of the face with the outward normal 
$\boldsymbol{n}_k$ for every $k=1,\dots,5$.
The set $\mathscr{L}(\boldsymbol{n}_1,\dots,\boldsymbol{n}_5)\subset\mathbb{R}^5$
can be referred to as a natural configuration space of convex 
polytopes (treated up to translations) with prescribed set 
$\{\boldsymbol{n}_1,\dots,\boldsymbol{n}_5\}$ 
of outward unit normals, when a polytope is determined by
the perimeters $L_1,\dots,L_5$ of its faces.

The main result of this article reads as follows:

\textbf{Theorem 5.} 
\textit{Let the vectors $\boldsymbol{n}_1,\dots,\boldsymbol{n}_5$ 
be given by the formulas} (1). \textit{Then the set 
$\mathscr{L}(\boldsymbol{n}_1,\dots,\boldsymbol{n}_5)\subset\mathbb{R}^5$ 
is not locally-analytic.}

From our point of view, Theorem 5 explains why a general 
existence theorem is not known which 
determines a convex polytope in $\mathbb{R}^3$ via unit 
normals and perimeters of its faces. 
The reason is that no analytic condition, similar to 
the condition (iii) in Theorem 4, does exist.

\textbf{2. Auxiliary constructions and preliminary results.}
Let $P\subset\mathbb{R}^3$ be a convex polytope such that 
the vectors $\boldsymbol{n}_1,\dots,\boldsymbol{n}_5$ 
defined by the formulas (1) (and no other vector) are the 
unit outward normals to the faces of $P$.
For $k=1,\dots,5$, denote by $\pi_k$ the 2-dimensional
plane in $\mathbb{R}^3$ containing the face of $P$ with the 
outward normal $\boldsymbol{n}_k$.

The straight lines $\pi_1\cap\pi_2$ and $\pi_1\cap\pi_3$
are parallel to the vector $\boldsymbol{e}_2=(0,1,0)$,
and the straight lines $\pi_1\cap\pi_4$ and $\pi_1\cap\pi_5$
are parallel to the vector $\boldsymbol{e}_1=(1,0,0)$.
Hence, the face $P\cap\pi_1$ is a rectangle.
Computing the angles between the vectors $\boldsymbol{n}_1$
and $\boldsymbol{n}_k$ for $k=2,\dots, 5$, we conclude that
the dihedral angle attached to any edge of the face $P\cap\pi_1$ 
is equal to $\pi/4$.
Now it is clear that the polytope $P$ can be of one of the three
types schematically shown in Figure~2. 
\begin{figure}
\begin{center}
\includegraphics[width=0.9\textwidth]{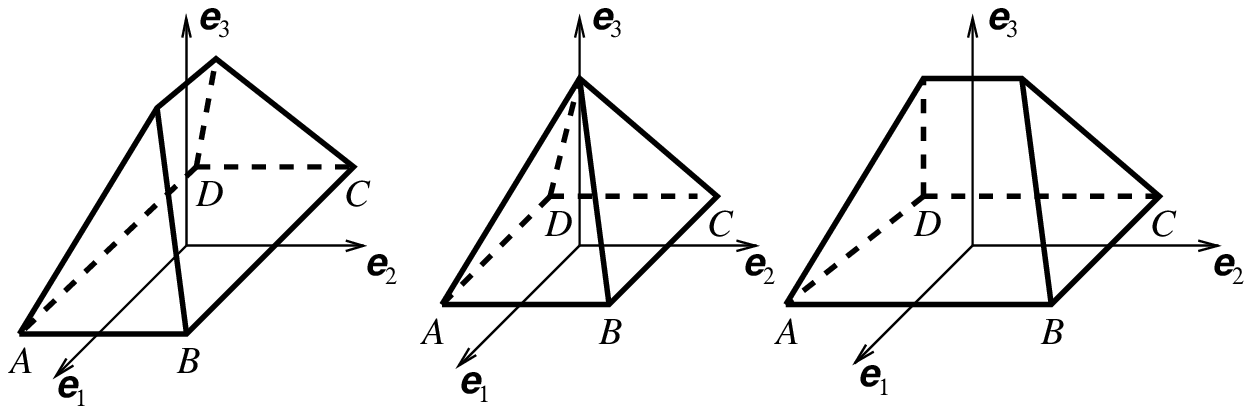}
\end{center}
\caption{Three types of polytopes with
outward normals $\boldsymbol{n}_1,\dots,\boldsymbol{n}_5$ 
defined by the formulas (1): Type I (left), Type II (center), 
and Type III (right)}
\end{figure}

In order to be more specific, we put by definition
$A=\pi_1\cap\pi_2\cap\pi_5$,
$B=\pi_1\cap\pi_2\cap\pi_4$,
$C=\pi_1\cap\pi_3\cap\pi_4$, and
$D=\pi_1\cap\pi_3\cap\pi_5$.
Denote by $2x$ the length of the straight line segment $AB$,
and by $2y$ the length of the straight line segment $BC$.
We say that the polytope $P$ is of Type I, if $x<y$;
is of Type II, if $x=y$; and is of Type III, if $x>y$.
Polytopes $P$ of Types I--III are shown schematically  in Figure~2. 

Denote by 
$\mathscr{L}_I(\boldsymbol{n}_1,\dots,\boldsymbol{n}_5)$ 
(respectively, by 
$\mathscr{L}_{II}(\boldsymbol{n}_1,\dots,\boldsymbol{n}_5)$
and $\mathscr{L}_{III}(\boldsymbol{n}_1,\dots,\boldsymbol{n}_5)$)
the set of all points $(L_1,\dots,L_5)\in \mathbb{R}^5$ such that
there exists a convex polytope $P\subset\mathbb{R}^3$ 
of Type I (respectively, of Type II or Type III) 
such that the vectors
$\boldsymbol{n}_1,\dots,\boldsymbol{n}_5$ 
(and no other vector) are outward unit normals
to the faces of $P$, and $L_k$ is the perimeter
of the face of $P$ with the outward normal $\boldsymbol{n}_k$
for every $k=1,\dots, 5$.
Below, we use also the following notation
\begin{eqnarray}
\boldsymbol{v}_I&=&(2, -(3+2\sqrt{3}), -(3+2\sqrt{3}), 5, 5),\nonumber\\
\boldsymbol{v}_{II}&=&(2(\sqrt{3}-1), 1,1,1,1),\nonumber\\
\boldsymbol{v}_{III}&=&(2,5,5,-(3+2\sqrt{3}), -(3+2\sqrt{3})).\nonumber
\end{eqnarray}
Note that in this paper the notation $AB$ can as well denote
the straight line segment and its length.

\textbf{Lemma 1.} 
\textit{Let the set $\{\boldsymbol{n}_1,\dots,\boldsymbol{n}_5\}$ 
of unit vectors in $\mathbb{R}^3$ be defined by the formulas} (1).
\textit{Then the following three statements are equivalent to each other:}

(i) $(L_1,\dots,L_5)\in\mathscr{L}_I(\boldsymbol{n}_1,\dots,\boldsymbol{n}_5)$;

(ii) $L_1=(2\sqrt{3}-3)L_2+L_4$, $L_2=L_3$, $L_4=L_5$, $L_4>L_2>0$;

(iii) $(L_1,\dots,L_5)=\alpha \boldsymbol{v}_I +\beta \boldsymbol{v}_{II}$
\textit{for some $\alpha, \beta\in\mathbb{R}$ 
such that} $\beta>(3+2\sqrt{3})\alpha>0.$

\textit{Proof}\,:
Suppose the statement (i) of Lemma 1 holds true.
In addition to the notation introduced above in Section~2, let 
$E=\pi_2\cap\pi_4\cap\pi_5$, $F=\pi_3\cap\pi_4\cap\pi_5$,
$G$ be the base of the perpendicular dropped from $E$ on the 
edge $AB$, 
$H$ be the base of the perpendicular dropped from $E$ on the face
$ABCD$, 
and $K$ be the base of the perpendicular dropped from $E$ on 
the edge $BC$, see Figure~3. 
\begin{figure}
\begin{center}
\includegraphics[width=0.45\textwidth]{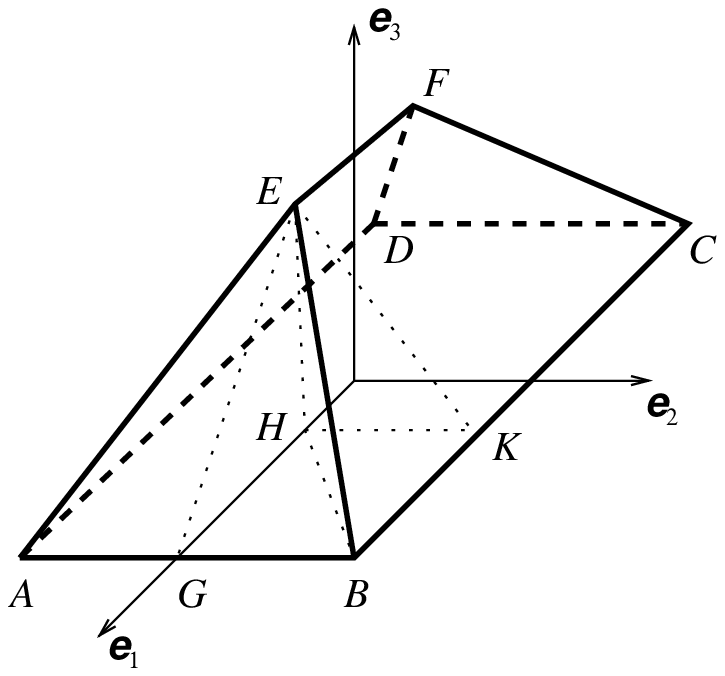}
\end{center}
\caption{
Elements of a polyhedron $P$ of Type I, needed for finding the 
relations between the perimeters $L_1,\dots, L_5$ of the faces}
\end{figure}
Using this notation, we obtain easily 
\begin{eqnarray}
BG&=&GH=EH=BK=x,\nonumber\\
BH&=&x\sqrt{2},\nonumber\\
BE&=&x\sqrt{3},\nonumber\\
EF&=&BC-2 BK=2y-2x,\nonumber\\
L_1&=&2AB+2BC=4x+4y,\\
L_2&=&L_3=AB+2BE=2(1+\sqrt{3})x,\\
L_4&=&L_5=BC+EF+2BE=2(\sqrt{3}-1)x+4y.
\end{eqnarray}
Eliminating $x$ and $y$ from the last three formulas, we get
$L_1=(2\sqrt{3}-3)L_2+L_4$.
Since the relations $L_2=L_3$, $L_4=L_5$, and $L_4>L_2>0$
are obvious for every polytope $P$ of Type I,
we conclude that the statement (i) implies the statement (ii).

Now suppose the statement (ii) of Lemma 1 holds true.
We, first, find the numbers $x$ and $y$ that satisfy the 
relations~(2)--(4) and, second, construct a convex 
polyhedron $P$ of Type I for which $L_1,\dots, L_5$ are 
the perimeters of the faces.

In accordance with (3), we put by definition 
$x=L_2(\sqrt{3}-1)/4$. Since $L_2>0$, $x>0$.
In accordance with (4), we put by definition
$y=L_2(\sqrt{3}-2)/4+L_4/4$.
Since $L_1=(2\sqrt{3}-3)L_2+L_4$, the condition (2) 
is satisfied: $4x+4y=(2\sqrt{3}-3)L_2+L_4=L_1.$
Using the inequalities $L_4>L_2>0$, we obtain
$4y=L_4+(\sqrt{3}-2)L_2>L_2+(\sqrt{3}-2)L_2=(\sqrt{3}-1)L_2>0$.
Hence, there exists a rectangle $ABCD$ in $\mathbb{R}^3$ 
such that the edge $AB$ is parallel to the vector 
$\boldsymbol{e}_2=(0,1,0)$ and its length is equal to $2x$,
and the edge $BC$ is parallel to the vector 
$\boldsymbol{e}_1=(1,0,0)$ and its length is equal to $2y$.

Denote by $\pi_1$ the plane that contains $ABCD$. 
Obviously, $\pi_1$ is perpendicular to $\boldsymbol{n}_1$.
Denote by $\pi_2$ the plane perpendicular to $\boldsymbol{n}_2$
and containing $AB$.
Denote by $\pi_3$ the plane perpendicular to $\boldsymbol{n}_3$
and containing $CD$.
Denote by $\pi_4$ the plane perpendicular to $\boldsymbol{n}_4$
and containing $BC$.
And denote by $\pi_5$ the plane perpendicular to $\boldsymbol{n}_5$
and containing $AD$.
The five planes $\pi_1,\dots,\pi_5$ determine a compact 
convex polyhedron with outward normals 
$\boldsymbol{n}_1,\dots,\boldsymbol{n}_5$.
Denote it by $P$. 
According to the statement (ii), $L_4>L_2$.
Hence, $4y=L_4+(\sqrt{3}-2)L_2>L_2+(\sqrt{3}-2)L_2=(\sqrt{3}-1)L_2=4x$.
Thus, $P\in\mathscr{L}_I$ and
the statement (ii) implies the statement (i).

So, we have proved that the statements (i) and (ii)
are equivalent to each other.

In order to prove that the statements (ii) and (iii)
are equivalent to each other, we observe that the equations 
$L_1=(2\sqrt{3}-3)L_2+L_4$, $L_2=L_3$, and $L_4=L_5$
from the statement (ii) define a 2-dimensional plane in 
$\mathbb{R}^5$. Denote this plane by $\lambda_I$.
Moreover, the vectors $\boldsymbol{v}_I$
and $\boldsymbol{v}_{II}$ constitute an orthogonal basis 
in $\lambda_I$.
This means that every vector $(L_1,\dots,L_5)\in\lambda_I$
can be uniquely written in the form 
$\alpha\boldsymbol{v}_I+\beta\boldsymbol{v}_{II}$.
Direct calculations show that the inequalities
$L_4>L_2>0$ from the statement (ii) are equivalent 
to the inequalities $\beta>(3+2\sqrt{3})\alpha>0$
from the statement (iii). \hfill $\square$

\textbf{Lemma 2.} 
\textit{Let the set $\{\boldsymbol{n}_1,\dots,\boldsymbol{n}_5\}$ 
of unit vectors in $\mathbb{R}^3$ be defined by the formulas} (1).
\textit{Then the following three statements are equivalent to each other:}

(i) $(L_1,\dots,L_5)\in\mathscr{L}_{II}(\boldsymbol{n}_1,\dots,\boldsymbol{n}_5)$;

(ii) $L_1=2(\sqrt{3}-1)L_2$, $L_2=L_3=L_4=L_5>0$;

(iii) $(L_1,\dots,L_5)=\gamma \boldsymbol{v}_{II}$
\textit{for some $\gamma\in\mathbb{R}$ 
such that} $\gamma>0.$

\textit{Proof} is left to the reader.
It can be obtained by arguments similar to those used above
in the proof of Lemma 1. 
But in fact, it is sufficient to observe that Lemma 2
is the limit case of Lemma 1 as $L_4$ approaches $L_2$.

\textbf{Lemma 3.} 
\textit{Let the set $\{\boldsymbol{n}_1,\dots,\boldsymbol{n}_5\}$ 
of unit vectors in $\mathbb{R}^3$ be defined by the formulas} (1).
\textit{Then the following three statements are equivalent to each other:}

(i) $(L_1,\dots,L_5)\in\mathscr{L}_{III}(\boldsymbol{n}_1,\dots,\boldsymbol{n}_5)$;

(ii) $L_1=L_2+(2\sqrt{3}-3)L_4$, $L_2=L_3$, $L_4=L_5$, $L_2>L_4>0$;

(iii) $(L_1,\dots,L_5)=\delta \boldsymbol{v}_{III} +\varepsilon \boldsymbol{v}_{II}$
\textit{for some $\delta, \varepsilon\in\mathbb{R}$ 
such that} $\varepsilon>(3+2\sqrt{3})\delta>0.$

\textit{Proof} is left to the reader.
It can be obtained by arguments similar to those used above
in the proof of Lemma 1. 
But in fact, it is sufficient to observe that if we rotate 
a polytope of Type III around the vector 
$\boldsymbol{e}_3=(0,0,1)$ to the angle $\pi/2$,
we get a polytope of Type I and can apply Lemma 1 to it. 

In the proof of Lemma 1, we denoted by $\lambda_I$ 
the 2-dimensional subspace in $\mathbb{R}^5$
which is spanned by the vectors 
$\boldsymbol{v}_I$ and $\boldsymbol{v}_{II}$.
Now we denote by $\lambda_{II}$ 
the 1-dimensional subspace in $\mathbb{R}^5$ 
spanned by $\boldsymbol{v}_{II}$
and denote by $\lambda_{III}$ 
the 2-dimensional subspace spanned by $\boldsymbol{v}_{II}$
and $\boldsymbol{v}_{III}$. 
  
\textbf{Lemma 4.} 
$\lambda_{II}=\lambda_I\cap\lambda_{III}$.

\textit{Proof}\,:
Each subspace $\lambda_I$, $\lambda_{II}$, and $\lambda_{III}$
contains $\boldsymbol{v}_{II}$.
Hence, $\dim(\lambda_I\cap\lambda_{III})\geqslant 1$. 
On the other hand, $\dim\lambda_I=\dim\lambda_{III}=2$. 
Hence, $\dim(\lambda_I\cap\lambda_{III})$ is equal to either 1 
or 2. 

Suppose $\dim(\lambda_I\cap\lambda_{III})=2$.
Then $\lambda_I=\lambda_{III}$.
Hence, the vectors $\boldsymbol{v}_I$, $\boldsymbol{v}_{II}$, 
and $\boldsymbol{v}_{III}$ are linearly dependant.
But this is not the case because the $3\times3$ minor 
composed of the first, third and fifth columns of the matrix
$$
\begin{pmatrix}
\boldsymbol{v}_I\\
\boldsymbol{v}_{II}\\
\boldsymbol{v}_{III}
\end{pmatrix}
=
\begin{pmatrix}
2&             -(3+2\sqrt{3})& -(3+2\sqrt{3})& 5&              5 \\
2(\sqrt{3}-1)& 1&              1&              1&              1\\
2&             5&              5&              -(3+2\sqrt{3})& -(3+2\sqrt{3})
\end{pmatrix}
$$
is non-zero. Hence, $\dim(\lambda_I\cap\lambda_{III})=1$, and
$\lambda_{II}=\lambda_I\cap\lambda_{III}$. \hfill $\square$

\textbf{3. Half-branches of analytic sets and the proof of Theorem 5.} 
Let $A$ be a one-dimensional analytic set, and $x\in A$. 
For every sufficiently small open ball $U$ with center $x$, 
$A\cap (U\diagdown\{x\})$ has a finite number of 
connected components $A_1,\dots,A_k$ such that $x$ belongs
to the closure of $A_j$ for every $j=1,\dots, k$.
These $A_j$ are called the half-branches of $A$ 
centered at $x$.
It is known that \textit{the number of half-branches 
of a one-dimensional analytic set centered at a point is even,} 
see, e.\,g. \cite{Su71}.

For completeness, we mention that a comprehensive exposition of a similar result for 
algebraic sets of dimension 1 may be found in \cite[Section~9.5]{BCR8}.

\textit{Proof of Theorem 5}\,:
Let $\Lambda$ be the straight line in $\mathbb{R}^5$
defined by the formula
$\Lambda=\{x\in\mathbb{R}^5\vert x=
\boldsymbol{v}_{II}+t\boldsymbol{v}_I \ 
\textrm{for some $t\in\mathbb{R}$}\}.$

Our proof is by contradiction.
Suppose the set 
$\mathscr{L}(\boldsymbol{n}_1,\dots,\boldsymbol{n}_5)\subset\mathbb{R}^5$ 
is locally-analytic.
Then  
$\Lambda\cap\mathscr{L}(\boldsymbol{n}_1,\dots,\boldsymbol{n}_5)$
is also locally-analytic.
Moreover, it is one-dimensional, contains the point
$\boldsymbol{v}_{II}$, and has only one half-branch centered at
$\boldsymbol{v}_{II}$. 

Let us explain the last statements in more details.
From the definition of polytopes of Types I--III we know that
$$
\mathscr{L}(\boldsymbol{n}_1,\dots,\boldsymbol{n}_5)=
\mathscr{L}_I(\boldsymbol{n}_1,\dots,\boldsymbol{n}_5)\cup
\mathscr{L}_{II}(\boldsymbol{n}_1,\dots,\boldsymbol{n}_5)\cup
\mathscr{L}_{III}(\boldsymbol{n}_1,\dots,\boldsymbol{n}_5).
$$
From Lemma 1 we know that 
$\mathscr{L}_I(\boldsymbol{n}_1,\dots,\boldsymbol{n}_5)$
is an angle on the 2-dimensional plane $\lambda_I\subset\mathbb{R}^5$.
From Lemma 3 we know that 
$\mathscr{L}_{III}(\boldsymbol{n}_1,\dots,\boldsymbol{n}_5)$
is an angle on the 2-dimensional plane $\lambda_{III}\subset\mathbb{R}^5$.
These angles are glued together along the ray
$\mathscr{L}_{II}(\boldsymbol{n}_1,\dots,\boldsymbol{n}_5)$ (see Lemma 2), and
no 2-dimensional plane contains the both of them (see Lemma 4). 
The line $\Lambda$ lies in the plane $\lambda_I$ and passes through the point
$\boldsymbol{v}_{II}$. Hence, for every sufficiently small open ball 
$U\subset\mathbb{R}^5$ with center $\boldsymbol{v}_{II}$, 
$$
U\cap\Lambda\cap\mathscr{L}(\boldsymbol{n}_1,\dots,\boldsymbol{n}_5)=
\{\boldsymbol{v}_{II}\}\cup(U\cap\Lambda\cap\mathscr{L}_I(\boldsymbol{n}_1,\dots,\boldsymbol{n}_5)).
$$
This formula means that we may obtain
$U\cap\Lambda\cap\mathscr{L}(\boldsymbol{n}_1,\dots,\boldsymbol{n}_5)$ in the following way:
first, we divide the straight line $\Lambda$ into two rays by the point
$\boldsymbol{v}_{II}$; then we observe that only one of these rays has at least one
common point with the angle $\mathscr{L}_I(\boldsymbol{n}_1,\dots,\boldsymbol{n}_5)$ 
and select that ray; at last, we intersect the ray selected with $U$. 

From this description, it is clear that 
$U\cap\Lambda\cap\mathscr{L}(\boldsymbol{n}_1,\dots,\boldsymbol{n}_5)$
is the half-branch of the locally-analytic set 
$\mathscr{L}(\boldsymbol{n}_1,\dots,\boldsymbol{n}_5)\subset\mathbb{R}^5$
centered at $\boldsymbol{v}_{II}$. 
Moreover, this is the only half-branch centered at $\boldsymbol{v}_{II}$. 
This contradicts to the fact that the number of half-branches 
of a locally-analytic set centered at a point is even, see \cite{Su71}.
\hfill $\square$

\textit{Remark}\,:
The proof of Theorem 5 provides us with a new, more technical, answer to the 
question of the title of this article.
As a part of the proof of Theorem 5, we demonstrated that the set 
$\mathscr{L}(\boldsymbol{n}_1,\dots,\boldsymbol{n}_5)\subset\mathbb{R}^5$ is not convex.
In Section 1, we mentioned that $\mathscr{L}(\boldsymbol{n}_1,\dots,\boldsymbol{n}_5)$
can be considered as a natural configuration space of convex 
polytopes (treated up to translations) with prescribed 
outward unit normals and perimeters of its faces.
The reader, familiar with the proof of Theorem 4 given in \cite[Chapter VII, \S 1]{Al05}, 
may remember that convexity of the analogous `natural configuration space'
$\mathscr{F}(\boldsymbol{n}_1,\dots,\boldsymbol{n}_m)\subset\mathbb{R}^m$
plays an important role in that proof.

\noindent{Victor Alexandrov}

\noindent\textit{Sobolev Institute of Mathematics}

\noindent\textit{Koptyug ave., 4}

\noindent\textit{Novosibirsk, 630090, Russia}

and

\noindent\textit{Department of Physics}

\noindent\textit{Novosibirsk State University}

\noindent\textit{Pirogov str., 2}

\noindent\textit{Novosibirsk, 630090, Russia}

\noindent\textit{e-mail: alex@math.nsc.ru}


\begin{thebibliography}{8}

\bibitem{Al37}{Alexandrov, A.D.}:
An elementary proof of the Minkowski and some other theorems on convex polyhedra (in Russian).
Izv. Akad. Nauk SSSR, Ser. Mat. No.~4, 597--606 (1937). 
JFM 63.1234.02

\bibitem{Al96}{Alexandrov, A.D.}:
Selected works. Part 1: Selected scientific papers.
Gordon and Breach Publishers, Amsterdam (1996).
MR1629804, Zbl 0960.01035

\bibitem{Al05}{Alexandrov, A.D.}:
Convex polyhedra. 
Springer, Berlin (2005).
MR2127379, Zbl 1067.52011

\bibitem{BCR8}{Bochnak, J.; Coste, M.; Roy, M.-F.}:
Real algebraic geometry.
Springer, Berlin (1998).
MR1659509, Zbl 0912.14023

\bibitem{Mi97}{Minkowski, H.}:
Allgemeine Lehrs{\"a}tze {\"u}ber die convexen Polyeder. 
G{\"o}tt. Nachr. 198--219 (1897).
JFM 28.0427.01

\bibitem{Mi11}{Minkowski, H.}:
Gesammelte Abhandlungen von Hermann Minkowski. Band I.
Teubner, Leipzig (1911).
JFM 42.0023.03

\bibitem{Pa08}{Panina, G.}:
A.D. Alexandrov's uniqueness theorem for convex polytopes and its refinements.
Beitr. Algebra Geom. \textbf{49}, No.~1, 59--70 (2008).
MR2410564, Zbl 1145.52007

\bibitem{Su71}{Sullivan, D.}:
Combinatorial invariants of analytic spaces. 
Proc. Liverpool Singularities-Sympos. I, 
Dept. Pure Math. Univ. Liverpool 1969--1970, 165--168 (1971).
MR0339241, Zbl 0227.32005

\end{thebibliography}
\end{document}